\documentclass[12pt]{article}
\usepackage{epsfig,makeidx,amsmath,amsfonts,latexsym,subfigure}
\usepackage{times,pifont}
\newtheorem{theorem}{Theorem}[section]
\newtheorem{prop}[theorem]{Proposition}
\newtheorem{lemma}[theorem]{Lemma}

\newcommand{\E}{{\mathbb E}}

\newcommand{\Var}{{\mathbb V}{\rm ar}}
\newcommand{\Cov}{{\mathbb C}{\rm ov}}
\newcommand{\RR}{{\mathbb R}}

\newcommand{\Z}{{\mathbb Z}}

\newcommand{\Pro}{{\mathbb P}}

\newcommand{\taum}{\tau_{{\rm mix}}}
\newcommand{\tauh}{\hat{\tau}}

\newcommand{\n}{\|}
\newcommand{\bi}{\begin{itemize}}
\newcommand{\ei}{\end{itemize}}
\newcommand{\be}{\begin{enumerate}}
\newcommand{\ee}{\end{enumerate}}
\newcommand{\beq}{\begin{equation}}
\newcommand{\eeq}{\end{equation}}
\newcommand{\beqa}{\begin{eqnarray*}}
\newcommand{\eeqa}{\end{eqnarray*}}
\newcommand{\ra}{\rightarrow}

\newcommand{\ep}{\epsilon}

\newcommand{\ov}{\overline}

\begin{document}

\title{Mixing times for the interchange process}

\author{Johan Jonasson
\thanks{Chalmers University of Technology and G\"oteborg University, Dept. of Mathematical Sciences, S-412 96 G\"oteborg, Sweden}
\thanks{email: jonasson@chalmers.se, url: http://www.math.chalmers.se/$\sim$jonasson}}

\date{May 2012}

\maketitle

\begin{abstract}
Consider the interchange process on a connected graph $G=(V,E)$ on $n$ vertices. I.e.\ shuffle a deck of cards by first placing one card at each vertex of $G$ in a fixed order and then at each tick of the clock, picking an edge uniformly at random and switching the two cards at the end vertices of the edge with probability $1/2$.
Well known special cases are the random transpositions shuffle, where $G$ is the complete graph, and the transposing neighbors shuffle, where $G$ is the $n$-path.
Other cases that have been studied are the $d$-dimensional grid, the hypercube, lollipop graphs and Erd\H os-R\'enyi random graphs above the threshold for connectedness.

In this paper the problem is studied for general $G$. Special attention is focused on trees, random trees and the giant component of critical and supercritical $G(N,p)$ random graphs. Upper and lower bounds on the mixing time are given. In many of the cases, we establish the exact order of the mixing time.
We also mention the cases when $G$ is the hypercube and when $G$ is a bounded-degree expander, giving upper and lower bounds on the mixing time.
\end{abstract}


\smallskip


\section{Introduction}

Card shuffling has, for several decades now, been one of the major playgrounds for
developing methods for the study of mixing rates for Markov chains.
One of the first card shuffling techniques to be studied was the random transpositions shuffle, where at each (discrete) time, two positions in the deck are chosen uniformly at random and the cards at those positions are switched with probability $1/2$.
(Strictly speaking, this is the lazy version of the random transpositions shuffle.)
If the deck consists of $n$ cards and starts out from a fixed permutation,
then an easy coupon collector's argument shows that the mixing time in total variation
is at least $n\log n$.
Diaconis and Shahshahani \cite{DS} proved that this is also sufficient,
not only in total variation, but also in $L^2$.
Matthews \cite{Matthews} later proved, via a strong uniform time argument, that
this is also sufficient in separation.
Hence the total variation norm exhibits a sharp threshold at time $n\log n$
and hence $\taum = (1+o(1))n\log n$.

Another shuffle that was studied in the early days, is the transposing neighbors shuffle, where at each time a position $i\in [n-1]$ in the deck is chosen uniformly and the card in that position is switched with the card in position $i+1$ with probability $1/2$.
Aldous \cite{Aldous} showed that the mixing time is at most of order $n^3\log n$ and at least of order $n^3$.
Much later Wilson \cite{Wilson1} introduced an important technique for lower bounding mixing time. This allowed him, among other things, to establish that
the true order of mixing for transposing neighbors is $n^3\log n$.
Wilson also considered neighbor transpositions on a grid and on the hypercube, proving that the order of mixing on the $\sqrt{n} \times \sqrt{n}$ grid is
$n^2\log n$ and that the order of the mixing time on $\Z_2^d$ is at least $n\log^2 n$, where
$n=2^d$.
Wilson conjectured that $n\log^2 n$ is the correct order of mixing on the hypercube. We partially confirm this by showing that mixing time is $O(n\log^3n)$.

\noindent {\bf Remark on notation.} Consider two expressions $f(n)$ and $g(n)$ as $n \ra \infty$. We will adopt the notation $f(n)=o(g(n))$ when $f(n)/g(n) \ra 0$. 
Writing $f(n) = O(g(n))$ means that that for some constant $C<\infty$, $f(n) \leq Cg(n)$ for all $n$. We write $f(n) = \Omega(g(n))$ for $g(n) = O(f(n))$ and
$f(n) =\omega(g(n))$ for $g(n) = o(f(n))$.  
We will also use the abbreviation "whp" for "with high probability". I.e.\ as a parameter, which is understood from the context, tends to $\infty$, an event that occurs with probability $1-o(1)$ is said to occur whp.

In the light of the above, it appears natural to consider neighbor transposition shuffling on a general connected $n$-vertex graph $G=(V,E)$, i.e.\ the card shuffle that at each time chooses uniformly at random an edge $e \in E$ and with probability $1/2$ switches the two cards at the end vertices of $e$.
Indeed, this general setup was considered by Aldous and Fill \cite{AF}, Chapter 14, but their analysis did not go further than what we have already mentioned, in terms of analyzing the process on particular $G$'s with respect to mixing time.
Aldous and Fill introduced the name {\em interchange process} for neighbor transposition shuffling on general graphs.
The interchange process has attracted a fair amount of interest in recent years. The most prominent result is perhaps the proof by Caputo, Liggett and Richthammer  \cite{CLR} of the
Aldous spectral gap conjecture (see e.g.\ David Aldous homepage), which states in our setting, that the spectral gap of the interchange process on $G$ equals the spectral gap of simple random walk on $G$ times $n$ for {\em all} $G$.
It is not clear if this has any direct relation to the mixing time results in this paper. As will be seen in our main result, Theorem \ref{tA} below, the relation between the mixing times of the interchange process and of random walk is usually not the same as for the spectral gap. Indeed, much of our work will go into proving that in many cases mixing time for the interchange process exceeds the mixing time (and the spectral gap) of simple random walk by a factor $n\log n$.
More on this will follow in the remark after Theorem \ref{tA}.

A recent contribution was given by Erikshed \cite{Erikshed}.
He considered the interchange process on lollipop graphs with $\Theta(n)$ vertices in the clique as well as in the handle and proved that
$\taum = \Theta(n^4\log n)$.
The following result establishes that this is the highest possible order of mixing for the interchange process on any graph.

\begin{prop} \label{pA}
Consider the interchange process on $G=(V,E)$. Let $m=|E|$ and $\rho$ be the radius of $G$. Then
\[\taum \leq (1+o(1))8m\rho n\log n.\]
\end{prop}

Erikshed also considered the interchange process on the Erd\H os-R\'enyi $G(n,p)$ random graph for $p > \log n/n$, the threshold for connectedness.
He proved that whp $\taum = \Omega(n\log n)$ and when $np=\Omega(n^\delta)$, $\delta>0$, then this is the true order of mixing.

In this paper, we consider shuffling on many other graphs. The following general result will be our main tool for lower bounding mixing time.

\begin{prop} \label{pB}
Assume that $1-\gamma$, $\gamma = o(1/n)$, is the second largest eigenvalue for the motion of a single card under the interchange process on $G$.
Let $\xi:V \ra \RR$ be a corresponding eigenvector with $\n\xi\n_2=1$ and assume that $\n\xi\n_1 \geq n^a$ for some $a>0$.
Then
\[\taum = \Omega(\gamma^{-1}\log n).\]
\end{prop}

The proof builds on Wilson's technique, introduced in \cite{Wilson1} (and developed further in \cite{Wilson2}, \cite{J1} and \cite{J2}).
As in Wilson's technique, we use the eigenvector to construct a test function, but since
we will not know the precise form of the eigenvector, we bound the variance of the test function in a different way.
The condition that the $L^1$-norm of the second eigenvector is not too low is a mild restriction. An eigenvector corresponding to the second largest eigenvalue of the single-card chain, is also known as a Fiedler vector of the graph Laplacian, $L(G)$ (in honor of the pioneering work of Fiedler, \cite{Fiedler73} and \cite{Fiedler75}) among graph theorists.
The Fiedler vector has turned out to be of significance in graph theory, most prominently in graph partitioning.
Intuitively, the $L^1$-norm of a Fiedler vector should be of order $\sqrt{n}$ and I conjecture that it is at least of order $n^a$ for all graphs.
However, even though many facts about the structure of Fiedler vectors are known, see e.g.\ \cite{BLS}, I have not been able to find a result of this kind.
Nevertheless, it is easy to establish the required lower bound on the $L^1$-norm in most of our applications;

\begin{lemma} \label{lC}
Let $G=(V,E)$ be an $n$-vertex connected graph. Suppose that $\xi:[n] \ra \RR$ is an eigenvector of $L(G)$ corresponding to an eigenvalue $\kappa=O(n^{-b})$ for
some $b>0$ (i.e.\ $\xi$ is an eigenvector for the random walk of a single card under the interchange process on $G$, corresponding to the eigenvalue $1-\kappa/2m$) and such that $\n\xi\n_2=1$.
Then $\n \xi \n_1 = \Omega(n^a)$ for some $a>0$.
\end{lemma}

Our main result is the following, giving (bounds for) the mixing time on several different graphs.
Here $\taum$ denotes the mixing time in terms of total variation and $\tauh$ is the convergence time in $L^2$.
The precise definitions will appear in Section 2.

\begin{theorem} \label{tA}
Consider the interchange process on $G=(V,E)$. Set $n:=|V|$.
\bi
\item[(a)] Let $d=\Theta(\log n)$ and let $G$ be the first $d$ generations of a rooted $r$-regular tree (i.e.\ when $V$ consists of all vertices within distance $d$ of a distinguished vertex of an $r$-regular tree), $3 \leq r = \Theta(1)$. Then
    \[\Omega(n^2\log n) =\taum \leq \tauh = O(n^2\log^2n).\]
    This also goes whp when $G$ is the first $d$ generations of a supercritical Galton-Watson tree, with finite variance offspring distribution, conditioned to survive.
\item[(b)] Let $G$ be a uniform random (labelled) tree on $n$ vertices. Then whp
    $\taum$ and $\tauh$ are both $\Theta(n^{5/2}\log n)$.
    This also goes whp for the $d=\Theta(n^{1/2})$ first generations of the incipient infinite tree of a critical Galton-Watson process whose offspring distribution has finite variance.
\item[(c)] Let $G$ be the giant component of a $G(N,(1+\ep)/N)$ graph, where $N=\Theta(n^{3/2})$ and $\ep=O(N^{-1/3})$. Then whp
    $\taum$ and $\tauh$ are both $\Theta(n^{5/2}\log n)$.
\item[(d)] Let $G$ be the giant
component of a $G(N,(1+\ep)/N)$ graph where $N=\ep^{-1}n$ and $\omega(N^{-1/3}) = \ep = o(1)$. Then whp
\[\taum=\Theta(n\ep^{-3}\log^2(\ep^2n)\log n).\]
The same thing holds for $\tauh$.
\item[(e)] Let $G$ be the giant component of a $G(N,(1+\ep)/N)$ graph where $N=\Theta(n)$ and $\ep=\Theta(1)$, Then whp
    $\taum$ and $\tauh$ are both $\Theta(n\log^3 n)$.
\item[(f)] When $G$ is a bounded degree expander, then $\Omega(n\log n) = \taum \leq \tauh = O(n\log^3n)$.
\item[(g)] When $G = Z_2^d$, $d=\log_2 n$, then
$\Omega(n\log^2n) = \taum \leq \tauh = O(n\log^3n)$.
\ei
\end{theorem}

{\bf Remark.} In (b)-(e), the mixing time of the interchange process on $G$ is of the order
$n\log n$ times the mixing time of random walk on $G$, (see
\cite{BKW}, \cite{DLP}, \cite{DHLP1}, \cite{FR} and \cite{NP}
for the mixing time of random walk on these graphs).
I believe this to be true for (a) as well.
On the other hand, in (f) the true mixing time appears to be $n\log n$ which is $\Theta(n)$ times the mixing time for random walk.
In (g), the mixing time for random walk is $d\log d$, so the ratio between the two mixing times is between $\Theta(n\log n/\log\log n)$ and $\Theta(n\log^2n/\log\log n)$.
It would be interesting to know what possible relationships one can have between the two mixing times.
It would also be interesting to compare mixing times for the interchange process with mxing times for the symmetric exclusion process on the same graphs, for which new results can be found in \cite{Oliveira}.
\smallskip

The rest of the paper is organized as follows:
the basic concepts are introduced in Section 2 and the proofs are given in Section 3.

\section{Preliminaries}

Let $S$ be a finite set and $\nu$ a signed measure on $S$ with $\nu(S)=0$
and let $\pi$ be a probability measure on $S$.
For $p \in [1,\infty)$, the $L^p$-norm of $\nu$ with respect to $\pi$ as
\[\n\nu\n_p = \Big( \sum_{i \in S} \Big|\frac{\nu(i)}{\pi(i)}\Big|^p \pi(i)\Big)^{1/p}.\]
The total variation norm of $\nu$ is given by
\[\n\nu\n_{TV} = \frac12 \sum_{i \in S}|\nu(i)| = \max_{A \subseteq S} \nu(A).\]
By Cauchy-Schwarz, $2\n\nu\n_{TV}=\n\nu\n_1 \leq \n\nu\n_2$.

When $\{X_t\}$ is an irreducible aperiodic Markov chain on $S$, with stationary distribution $\pi$, we define the mixing time as
\[\taum = \min\Big\{t:\n \Pro(X_t \in \cdot)-\pi\n_{TV} \leq \frac14\Big\}.\]
The convergence time in $L^2$ is given by
\[\tauh = \min\Big\{t:\n\Pro(X_t \in \cdot) - \pi\n_2 \leq \frac12\Big\}.\]
It follows that $\taum \leq \tauh$.

All upper bounds in this paper rely on the comparison technique of Diaconis and Saloff-Coste \cite{DS}.
Let $\{X_t\}$ and $\{Y_t\}$ be two random walks on $S_n$ generated by the symmetric probability measures $\mu$ and $\mu_0$ respectively.
Let $E$ and $E_0$ be the supports of $\mu$ and $\mu_0$ respectively.
For each $y \in E_0$, choose a representation $y=x_1x_2\ldots x_k$, $x_i \in E$ and $k$ odd, and write $|y|=k$.
For each $x \in E$ and $y \in E_0$, let $N(x,y)$ be the number of times that $x$ appears in the chosen representation of $y$.
Define
\[A_* = \max_{x \in E}\frac{1}{\mu(x)} \sum_{y \in E_0}|y|N(x,y)\mu_0(y).\]
Then the following is a special case of Lemma 5 of \cite{DS}
\begin{lemma} \label{lA}
The distance in $L^2$ from stationarity of $\{X_t\}$ and $\{Y_t\}$ respectively relate in the following way.
\[\n\Pro(X_t \in \cdot)-\pi\n_2^2 \leq n!e^{-\lfloor t/A_* \rfloor}+
\n \Pro(Y_{\lfloor t/A_* \rfloor} \in \cdot)-\pi\n_2^2.\]
\end{lemma}
The most important consequence of the lemma is that, provided that
$\tauh_0 = \Omega(n\log n)$,
\[\tauh = O(A_* \tauh_0),\]
where $\tauh$ and $\tauh_0$ are the convergence times in $L^2$ of the two chains.
Of course, one also needs a chain $\{Y_t\}$ to compare with for which $\tauh_0$ is known. When it comes to card shuffling, the by far most common chain to compare with is the random transpositions shuffle. The following sharp result is due to Diaconis and Shahshahani \cite{DSh}.

\begin{lemma} \label{lB}
Let $\{Y_t\}$ be the random transposition shuffle. There exists a constant $C<\infty$ such that for $t=\lfloor n(\log n+c)\rfloor$,
\[\n\Pro(Y_t \in \cdot)-\pi\n_2^2 \leq Ce^{-2c}.\]
\end{lemma}

Before starting with the proofs, recall the extremal characterization of eigenvalues. For our purposes it suffices to recall the special case when $A$ is the symmetric transition matrix of an irreducible aperiodic Markov chain on the finite state space $S$. Then with $\kappa_2$ denoting the second largest eigenvalue of $A$, we have for the spectral gap $\gamma:=1-\kappa_2$
\[\gamma = \min\Big\{\xi^T(I-A)\xi: \sum_{v \in S}\xi(v)=0,\n \xi\n_2=1\Big\}\]
and a vector $\xi$ for which the minimum is attained is eigenvector for $\kappa_2$.
In this paper this will be applied on the case where $A$ is the transition
matrix for the motion of a single card under neighbor transpositions on $G=(V,E)$.
The neighbor transpositions shuffle on $G$ is the random walk on the symmetric group $S_n$ of permutations of $n$ cards, whose state changes from time $t$ to time $t+1$ in the following way. 
An edge $e \in E$ is chosen uniformly at random, independently of choices made at other times. Then, by another independent choice, the state at time $t+1$ is the composition of the state at time $t$ with $(i\,\,j)$, where $i$ and $j$ are the cards at the end vertices of $e$ at time $t$,
with probability $1/2$ and identical with the state at time $t$ with probability $1/2$.
This entails that any given single card makes a random walk on $V$ with the behavior that given that the card is at vertex $u$ at time $t$, it will at time $t+1$ be at $u$ with probability $1-d_u/(2|E|)$ and at $v$ with probability $1/(2|E|)$ for each neighbor $v$ of $u$.
Hence, when $A$ is the transition matrix for the single card random walk, we have
\begin{eqnarray*}
\xi^T(I-A)\xi & = & \frac{1}{2m}\sum_v\xi(v)\Big(d_v\xi(v)-\sum_{u \sim v}\xi(u)\Big) \\
& = & \frac{1}{2m}\sum_{\{u,v\} \in E}(\xi(u)-\xi(v))^2
\end{eqnarray*}
where $d_v$ is the degree of vertex $v$ and $u \sim v$ denotes that $\{u,v\} \in E$.

\section{Proofs}

\noindent {\bf Proof of Proposition \ref{pA}.}
This is a straightforward application of Lemma \ref{lA} and Lemma \ref{lB}.
Let $\{X_t\}$ be the interchange process on $G$ and let $\{Y_t\}$ be the random transpositions shuffle.
Let $y=(u\,\,v) \in E_0$ be an arbitrary transposition.
Let $w_0w_1\ldots w_r$ be a shortest path between $u$ and $v$ (where $w_0=v$ and $w_r=v$).
Then
\[(u\,\,v)=(w_0 \,\,w_1)(w_1\,\,w_2)\ldots (w_{r-2}\,\,w_{r-1})
(w_{r-1}\,\,w_r)(w_{r-2}\,\,w_{r-1})\ldots (w_0\,\,w_1).\]
Hence $|y| \leq 2r-1 < 2\rho$ and $N(x,y) \leq 2$.
Since $\mu(x) = 1/(2m)$ for any $x \in E$, $x \neq id$, Lemma \ref{lA} entails that $\tauh \leq (1+o(1))8m\rho\tauh_0)$, which by Lemma \ref{lB} is $(1+o(1))8m\rho n\log n)$.
\hfill $\Box$

\bigskip

\noindent {\bf Proof of Proposition \ref{pB}.}
Let $X_0$ be a fixed starting state and
let $\xi:V \ra \RR$ be an eigenvector corresponding to the eigenvalue $1-\gamma$ with $\sum_v\xi(v)^2=1$.
This means that if $X_t$ is the state of the deck and $X^i_t$ is the position of card $i$ at time $t$, then
\[\E[\xi(X^i_{t+1})|X_t]=(1-\gamma)\xi(X^i_t).\]
Let $\phi(X_t)=\sum_i\xi(X^i_t)$ where the sum is taken over the $i$'s for which
$\xi(X^i_0)>0$.
Then $\phi$ is an eigenvector for (the transition matrix of) $\{X_t\}$.
By assumption $\phi(X_0) \geq n^a$.
We claim that the terms are negatively correlated.
Indeed, writing $\Delta_i=\xi(X^i_t)-\xi(X^i_{t-1})$, we have
\[\E[\Delta_i\Delta_j|X_{t-1}] \leq 0\]
since either $\Delta_i$ or $\Delta_j$ is 0 unless card $i$ and card $j$ are adjacent in $X_t$ and the edge between them gets chosen for the update, in which case $\Delta_i\Delta_j \leq 0$.
Note also that $\E[\Delta_i|X_{t-1}] = -\gamma\xi(X^i_{t-1})$.
Hence
\beqa
\E[\xi(X^i_t)\xi(X^j_t)] &=& \E[\E[(\xi(X^i_{t-1})+\Delta_i)(\xi(X^j_{t-1})+\Delta_j)|X_{t-1}]] \\
&\leq& (1-2\gamma)\E[\xi(X^i_{t-1})\xi(X^j_{t-1})].
\eeqa
By induction,
\[\E[\xi(X^i_t)\xi(X^j_t)] \leq (1-2\gamma)^t\xi(X^i_0)\xi(X^j_0).\]
Similarly
\[\E[\xi(X^i_t)]\E[\xi(X^j_t)] = (1-\gamma)^{2t}\xi(X^i_0)\xi(X^j_0)
\geq (1-2\gamma)^t\xi(X^i_0)\xi(X^j_0).\]
Hence $\Cov(\xi(X^i_t),\xi(X^j_t)) \leq 0$.
Using this we have that
\[\Var(\phi(X_t)) \leq \sum_i\E[\xi(X^i_t)^2] < 1.\]
Since, by assumption, $\E[\phi(X_t)] = (1-\gamma)^t\phi(X_0) = \Omega(n^{a-b})$ for
$t:=b\gamma^{-1}\log n$, $0<b<a$, Chebyshev's inequality entails that
\[\Pro(\phi(X_t) < \frac12 n^{a-b}) = o(1).\]
On the other hand, if $X_\infty$ denotes a deck at uniformity, we have
$\Var(\phi(X_\infty))<1$ and $\E[\phi(X_\infty)]=0$ and hence
\[\Pro(\phi(X_\infty) > \frac12 n^{a-b})=o(1).\]
\hfill $\Box$

\bigskip

\noindent {\bf Proof of Lemma \ref{lC}.}
Assume for contradiction that $\n\xi\n_1 =O(n^a)$ for all $a>0$.
Then $\max_v \xi(v) = \Omega(n^{-a})$ for all $a>0$ (assuming that we have chosen $\xi$ so that $\max_v\xi(v) \geq -\min_v\xi(v)$).
By assumption and the extremal characterization of eigenvalues
\[O(n^{-b}) = \sum_{u \sim v}(\xi(u)-\xi(v))^2\]
so in particular $\max_{u \sim v}|\xi(u)-\xi(v)| = O(n^{-b/2})$.
Then, if $v_0$ is chosen so that $\xi(v_0) = \Omega(n^{-b/8})$, all vertices $v$ within distance $n^{3b/8}/2$ of $v_0$ have $\xi(v) \geq n^{-b/8}/2$.
If the diameter of $G$ is less than $n^{3b/8}/2$, this implies that $\xi(v)>0$ for all $\xi$, a contradiction since $\sum_v\xi(v)$ must be $0$.
If the diameter is at least $n^{3b/8}/2$, then we have $\n\xi\n_1 \geq n^{3b/8}/2 \cdot n^{-b/8}/2 = n^{b/4}/4$, a contradiction. \hfill $\Box$

\bigskip

\noindent {\bf Proof of Theorem \ref{tA}(a).}
The upper bounds follow immediately from Proposition \ref{pA}.

For the lower bounds, consider first the binary tree.
Define $\eta:V \ra \RR$ by $\eta(v)=0$ when $v$ is the root, $\eta(v)=1-2^{-k}$
when $v$ is at distance $k$ from the root, in the left third of the tree,
$\eta(v)=-(1-2^{-k})$ for vertices at distance $k$ from the root, in the right third of the tree and $\eta \equiv 0$ on the middle third.
Expressed differently, $\eta(v)$ is the electric potential at $v$ when we apply the potential $+1$ to the leafs of the left third of the tree and $-1$ to the vertices of the right third and regard all edges as unit resistors.

Then $\sum_v\eta(v)=0$ and $\sum_v\eta(v)^2 = \Theta(n)$.
Also, since $\eta$ is harmonic off the leafs, $d_v\eta(v) - \sum_{u \sim v}\eta(u)=0$ when $v$ is not a leaf.
For leafs $v$, we have
\[\eta(v)\Big(d_v\eta(v) - \sum_{u \sim v}\eta(u)\Big) < \frac{4}{2^d} < \frac{8}{n}.\]
Hence
\[\gamma \leq \frac{n\Theta(1/n)}{m\Theta(n)} = \Theta(1/n^2).\]
Plugging this into Proposition \ref{pB} and using Lemma \ref{lC} gives
\[\taum = \Omega(n^2\log n)\]
as desired.

We omit the analogous proof for $r>3$, so let now $G$ be a supercritical Galton-Watson tree conditioned on survival, where
$\mu^d = \Theta(n)$ and $\mu>1$ is the expectation of the offspring distribution.
Let $u$ be a vertex with at least two children whose progeny survives, in the first generation where such a vertex exists.
Let $u_1$ and $u_2$ be two arbitrarily chosen such children of $u$ and let $A_1$ and $A_2$ be the sets leafs in the progeny of $u_1$ and $u_2$ respectively.
Let $\eta(v)$, $v \in V$, be the electric potential in $v$, with boundary condition $\eta \equiv +1$ on $A_1$ and $\eta \equiv -1$ on $A_2$.
Using well-known facts about Galton-Watson processes, we have whp that,
$|A_1| = \Theta(\mu^d)=\Theta(n)$, $|A_2|=\Theta(n)$ and that the effective resistance between $A_1$ and $A_2$ is $\Theta(1)$.
Hence the total current flowing through edges incident to vertices of $A_1$ and $A_2$, is $\Theta(1)$, i.e.\ the sum of the absolute potential differences over these edges is $\Theta(1)$.
Since $\eta$ is harmonic off $A_1 \cup A_2$, we get
\[\sum_v\eta(v)\Big(d_v\eta(v) - \sum_{u \sim v}\eta(u)\Big) = \Theta(1).\]
We cannot, however, plug this into Proposition \ref{pB}, since typically $\sum_v\eta(v)$ is not $0$.
However replacing $\eta$ with $\phi := \eta-\overline{\eta}$, does not change the expression. Also, since $|A_i|=\Theta(n)$, it is clear that $\sum_v\phi(v)^2 = \Theta(n)$.
As above, this gives $\gamma = \Theta(1/n^2)$ and hence, by Proposition
\ref{pB} and Lemma \ref{lC},
\[\taum = \Omega(n^2\log n).\]
\hfill $\Box$

\bigskip

\noindent {\bf Proof of Theorem \ref{tA}(b).}
Let first $G=(V,E)$ be the first $d$ generations of the incipient infinite tree of a critical Galton-Watson process.
It is well-known that whp, $G$ contains $\Theta(d^2)$ vertices, i.e.\
$d=\Theta(n^{1/2})$.

The upper bound on $\taum$ now follows from Proposition \ref{pA}.

Let $A$ be the set of leafs at distance $d$ from the root. By \cite{FK}, Proposition 1.1 and \cite{CK}, Lemma 3.1, it follows that
whp the effective resistance between the root and $A$ is $\Theta(n^{1/2})$.
Let $u$ be a vertex closest to the root for which at least two children are connected to $A$. Pick arbitrarily two such children, $u_1$ and $u_2$ and let
$A_i$ be the set of vertices of $A$ which are in the progeny of $u_i$, $i=1,2$.
Whp $|A_i|=\Theta(d)$, see e.g.\ \cite{CK}, Proposition 2.6, (since the finite
variance condition on the offspring distribution makes sure that
the survival probability up to generation $d$ of the unconditioned Galton-Watson process is $\Theta(1/d)$), and the effective resistance between $A_1$ and $A_2$ is at least $\Theta(d)$. Apply the potential $+1$ to $A_1$ and $-1$ to $A_2$ and
let $\eta(v)$ be the potential at $v$, $v \in V$.
Letting $\phi=\eta-\ov{\eta}$, we get $\sum_v\phi(v)^2=\Theta(n)$ and
\[\sum_v\Big(d_v\phi(v)-\sum_{u \sim v}\phi(u)\Big) = \Theta(1/n^{1/2})\]
since the total current flowing into $A_2$ is $\Theta(n^{-1/2})$.
Hence $\gamma = \Theta(1/n^{5/2})$, so by Proposition \ref{pB},
\[\taum=\Omega(n^{5/2}\log n).\]

Now we turn to the random tree, $G$, on $n$ vertices.
It is well-known that the diameter of $G$ is whp $\Theta(n^{1/2})$,
so the upper bound again follows from Proposition \ref{pA}.

By \cite{AP}, Proposition 1, the law of $G$ is identical to the law of the Galton-Watson tree with Poisson(1) offspring distribution, conditioned to contain exactly $n$ vertices.
Let $d=Cn^{1/2}$. Whp $G$ contains $\Theta(n^{1/2})$ vertices at distance $d$ from the root for small enough $C$. Fix $u$, $u_1$ and $u_2$ as above, such that $u$ is a vertex closest
to the root for which at least two children are connected to level $d$.
Continue to mimic the above argument, letting $A_i$ be the set of vertices at level $d$ in the progeny of $u_i$.
Define $\eta$ and $\phi$ as above. Then $\sum_v\phi(v)^2 = \Theta(n)$ whp and by
\cite{LW},
$G$ is stochastically dominated by the incipient infinite cluster of the same Galton-Watson process.
(Indeed, it is shown in \cite{LW} that the trees conditioned on having exactly $n$ vertices, $n=1,2,\ldots$, is stochastically increasing in $n$.
This is shown for binomial$(k,1/k)$ offspring distribution for any $k$, so the Poisson(1) offspring distribution follows as a limiting case.)
Hence it follows that the effective resistance between $A_1$ and $A_2$ is whp $\Theta(n^{1/2})$.
The proof is now completed by an appeal to Proposition \ref{pA}.
\hfill $\Box$

\bigskip

\noindent {\bf Proof of Theorem \ref{tA} (c).}
Let $G=(V,E)$ be the giant component of the Erd\H os-R\'enyi random graph $G(N,(1+\ep)/N)$, $\ep = O(1/N^{1/3})$.
It is a classical result that whp $|V|=\Theta(N^{2/3})$. It is also known that $G$ is a tree, uniformly distributed among all trees on $V$, with a probability bounded away from $0$ and $1$.
It is easy to see that if $G$ is not a tree, then it contains whp only $O(1)$
cycles.
(For detailed results on the giant component and random walks on the giant component for these parameters, see \cite{NP}.)
Hence the results follow from the proof of (b) with only very minor modifications.
\hfill $\Box$

\bigskip

Before the proof of (d) and (e), we introduce a random graph model due to Ding et el.\ \cite{DHLP1} (for $\ep = o(1/N^{1/4}$)) and \cite{DHLP2} (extended version for $\ep=o(1)$) and further explored in \cite{DLP}.
The model is {\em contiguous} to the giant component of the $G(N,(1+\ep)/N)$ model,
$\omega(1/N^{1/3}) = \ep = o(1)$, in the sense that any whp property for one is also whp for the other.
The model is given by a three-step procedure as follows.
Here $\mu$ is the conjugate of $1+\ep$, i.e.\ $\mu<1$ and $\mu e^{-\mu}=(1+\ep)e^{-(1+\ep)}$.

\bi
\item[(1)] Let $\Gamma$ be a normal random variable with mean $1+\ep-\mu$ and variance $1/(\ep n)$. Let $D_u$, $u \in [N]$, be iid Poisson($\Gamma$) random
    variables, conditioned that $\sum_{u: D_u \geq 3}D_u$ is even.
    Let $M_k$, $k=1,2,\ldots$, be the number of $u$ such that $D_u=k$ and
    $M=\sum_{k \geq 3}M_k$.
    Select a random multigraph, $H$, on $M$ vertices uniformly among all multigraphs with $M_k$ vertices of degree $k$, $k \geq 3$.

\item[(2)] Let $K$ be the multigraph obtained from $H$ by replacing each edge $e \in E(H)$ with a path of independent geometric($1-\mu$) length.

\item[(3)] Attach to each vertex $v \in V(K)$ an independent Galton-Watson tree with Poisson($\mu$) offspring distribution (with $v$ as its root).
\ei

Note that the variance of the normal distribution in (1) is of smaller order than the square of the mean, so whp $\Gamma=\Theta(\ep)$.
This entails that $H$ contains whp $\Theta(\ep^3 N)$ vertices and edges and hence
$K$ has whp $\ep^2N$ vertices.

\bigskip

\noindent {\bf Proof of Theorem \ref{tA}(d)}.
This proof is fairly long, so we split into four parts: lower bound (general), lower bound for $\ep^{-1}$ of polynomial order, lower bound for $\ep^{-1}$ of subpolynomial order and upper bound.

\smallskip

\noindent {\bf Lower bound.}
Write $F$ for the graph constructed from (1)-(3) above.
We start with the lower bound.
Write $\mu=1-\delta$ and note that $\delta=(1+o(1))\ep$.
Let also
\[k_0=\ep^{-1}\log(\ep^3N/a)\]
where $a$ is a constant independent of $N$ and $\ep$.
Consider a Poisson($\mu$) tree. Let $Z_r$ be the size of the $r$'th generation and $T$ the total size of such a tree.
Condition now on that such a tree survives for at least $k$ generations, by conditioning on the {\em backbone} of the tree, i.e.\ the leftmost path of length $k$,  $P=v_0v_1 \ldots v_k$ from the root to generation $k$.
(To properly make sense of the word "leftmost", associate with each vertex, $v$, of a Poisson($\mu$) tree, an independent Poisson process of intensity $\mu$ on the unit interval and let $v$ have one child per point of this process. Then order the children from left to right according to the positions of their respective Poisson points.)
Then from each $v_i$ on the backbone, an independent tree emanates, having Poisson($\mu$) offspring in all but the first generation, for which the offspring is Poisson($\mu(1-X_i)$), where $X_i$ is the position in the unit interval of the Poisson point associated with $v_i$.
Let  $Z^{(i)}_r$ be the size of generation $r$ and $T^{(i)}$ the total size of the tree emanating from $v_i$, including $v_i$ so that $T=\sum_{r=0}^{k}T^{(i)}$.
Then
\[\E[Z^{(i)}_r|P] = (1-X_i)\mu^r\]
and
\[\E[T^{(i)}|P]= (1-X_i)\sum_{r=0}^\infty\mu^r=\frac{1-X_i}{\delta}.\]
Since $P$ is whp such that $X_i$ is less than, say, $2/3$ for more than half the $i$'s, we get whp
\beq \label{eb}
\E[Z_k|P] = \sum_{i=0}^k \E[Z^{(i)}_{k-i}|P] = \Theta(1/\delta)
\eeq
and hence $\E[T|P] = \Theta(k/\delta)$.
By (\ref{eb})
\[\Pro(Z_k>0) = \frac{\E[Z_k]}{\E[Z_k|Z_k>0]} = \frac{\Theta((1-\delta)^k)}{\Theta(1/\delta)} =
\Theta\Big(\delta(1-\delta)^k\Big)\]
which is at least $C/(\ep^2N)$ for any large $C$ if $k=k_0$ and the constant $a$ in the definition of $k_0$ is small enough.
Since there are whp $\Theta(\ep^2 N)$ vertices in $K$, at least two of the trees, attached to $K$ in step (3) above, will whp survive for $k_0$ generations.

Also, by the standard formula for the variance of the generation sizes of a Galton-Watson tree (modified for the first generation), whp,
\[\Var(Z^{(i)}_r|P) \leq \frac{\mu^r}{\delta}.\]
Using that
\[\Cov(Z^{(i)}_r,Z^{(i)}_s|P) \leq
(\Var(Z^{(i)}_r|P) \Var(Z^{(i)}_s|P))^{1/2} \leq \frac{\mu^{r/2}\mu^{s/2}}{\delta},\]
we get
\[\Var(T^{(i)}|P) \leq \frac{1}{\delta^2}\]
so that
\[\Var(T|P) \leq \frac{k_0}{\delta^2}.\]
Hence by Chebyshev's inequality, we have whp
\[T=\Theta(k_0/\delta) = \Theta(\ep^{-2}\log(\ep^3N)).\]

Now pick uniformly at random, two attached trees that survive up to generation $k_0$
and fix in each of them a leaf at distance at least $k_0$ from the root of it's tree.
Apply the potential $+1$ to one of these leafs and $-1$ to the other.
Let $\eta(v)$, $v \in V(F)$, be the resulting potential at $v$ and,
again, $\phi=\eta-\ov{\eta}$.
By the above considerations, $\sum_v\phi(v)^2 = \Omega(\ep^{-2}\log(\ep^3N))$.
We also get that
\[\sum_v\eta(v)\Big(d_v\eta(v) - \sum_{u \sim v}\eta(u)\Big) = \Theta(1/k_0).\]
Plugging into the extremal characterization of the second eigenvalue gives, for the random walk of a card,
\[\gamma^{-1} = \Omega\Big(n\ep^{-3}\log^2(\ep^2 n)\Big)\]
since whp $n=\ep N$.

\smallskip

\noindent {\bf Lower bound for $\ep^{-1}$ of polynomial order.}
Since we here assume that $\ep^{-1} = \Omega(n^a)$ for some $a>0$, we appeal to Lemma \ref{lC}.
By Proposition \ref{pB}, it then follows that
\[\taum = \Omega\Big(n\ep^{-3}\log^2(\ep^2n)\log n\Big),\]
establishing the desired result.

\smallskip

\noindent {\bf Lower bound for $\ep^{-1}$ of subpolynomial order.}
Here more work is needed, since unfortunately, I have not been able to establish the required condition on the $L^1$-norm of the second eigenvector in Proposition \ref{pB}.
Nevertheless, a lower bound of $\Theta(n\ep^{-3}\log^2(\ep^2n)\log n) = \Theta(n\ep^{-3}\log^3n)$ can be established, using other arguments.
Note that by the above construction, $F$ contains $\Theta(n\ep) = \Omega(n^b)$, $b<1$, trees attached in (3).
By the above arguments, whp, at least, say $n^{2/3}$ of these trees will survive up to generation $\alpha\ep^{-1}\log n$ and contain at least $\alpha\ep^{-2}\log n$ vertices,
provided that $\alpha$ is taken sufficiently small.
For each such tree, $T_i$, $i=1,\ldots,n^{2/3}$, pick a vertex $h=h_i$ half way (or half way plus $1/2$) between the root, $r=r_i$, of the tree and a leaf, $l=l_i$, most distant from $r$.
Consider the random walk performed by the card, $c=c_i$, that starts in $h$.
For vertices $u \in T_i$, let $\tau_u$ be the first time that $c$ visits $u$.
By the correspondence between random walks and electrical networks, $\Pro(\tau_l < \tau_r) = 1/2$.
Hence, the probability that $c$ walks to $l$ and back to $h$ at least $k$ times before it visits $r$, is $2^{-k}>2n^{-1/3}$ for $k=(1/3)\log n$. Write $A=A_i$ for this event.
Note also for later, that the probability that $c$ visits $l$ more than $2\log n$ times before it visits $r$ is $o(1/n)$, so whp this will not happen for any card in the deck.

Conditioned on $A$, each of the $k$ return trips to $l$ takes an independent time, which is the sum of the time taken to reach $l$ and the time to go from $l$ back to $h$.
The conditional expectation of the first of these terms is affected by the fact that we condition on hitting $l$ before $r$.
The distribution of the second term however, is unaffected by conditioning on $A$, since a random walk from $l$ must necessarily visit $h$ before it visits $r$.
Hence the expected time of a return trip is bounded from below by $\E_l\tau_h$ (in standard notation).
By the hitting time formula from \cite{HF}, the tree version of Tetali's \cite{Tetali} general hitting time formula, an exact expression for $\E_l\tau_h$ is
\[\E_l\tau_h = m\Big(d(l,h) + 2\sum_{x \in P}|E(G_x)|d(x,h)\Big).\]
Here $d(u,v)$ is the distance between $u$ and $v$, $P$ is the set of vertices of the unique path between $l$ and $h$ and $G_x$ is component containing $x$ of the graph obtained from $F$ by removing the neighbors in $P$ of $x$.
From the above arguments, it follows that whp a non-vanishing fraction of the vertices in $x \in P$ have $|E(G_x)| \geq \beta \ep^{-1}$ (for sufficiently small $\beta$) and are at distance $\geq (\alpha/3)\ep^{-1}\log n$ from $h$. Hence, whp,
\[\E_l\tau_h=\Theta(m\ep^{-3}\log^2 n).\]
This entails that $\E_l[\tau_r|A] = \Theta(n\ep^{-3}\log^3n)$ (since $m=\Theta(n)$ whp).
We also have
\[\Var_h(\tau_r|A_i) = \frac13\log _n \Var(C(h,l)|\tau_l<\tau_r)) + \Var_{v_i}(\tau_r)\]
where $C(u,v)$ is the commute time between $u$ and $v$, i.e.\ the time taken for a return trip from $u$ to $v$ and back.
Let $h^*$ be the vertex which maximizes $\E C(x,l)$ over $x \in T_i$. Here it is important to point out that here and in what follows until the order of $\Var_h(\tau_r|A)$ is established, we consider the random walk {\em restricted to $T_i$}.
Let $b = \E C(h^*,l)/ \E C(h,l)$ which is of constant order, since the effective resistance, i.e.\ distance, between $h$ and $l$ is of maximal order in $T_i$.
We get
\[\Pro(C(h,l) \geq kb\E C(h,l)|\tau_l<\tau_r) \leq 2\Pro(C(h,l) \geq k\E C(h^*,l)) \leq 2 \cdot 2^{-k/2}\]
since each round of $2\E C(h^*,l)$ has, by Markov's inequality, a probability of at least $1/2$ to finish the return trip.
Hence
\[\Var_h(C(h,l)|\tau_l<\tau_r) \leq \E_h[C(v_i,l)^2|\tau_l<\tau_r] = O(n^2\ep^{-6}\log^6 n).\]
In the same way, this also goes for the term $\Var_h(\tau_r)$, so summing up gives
\[\Var_{v_i}(\tau_r|A_i) = O(n^2\ep^{-6}\log^7 n).\]
By the above consideration, whp, $A_i$ will occur for at least $n^{1/3}$ of the $T_i$'s. Let $S$ be the sum of the $\tau_{r_i}$'s from the trees $T_i$ for which $A_i$ occurs.
Then (whp)
\[\E S = \Theta(n^{4/3}\ep^{-3}\log^3 n)\]
and, since the $A_i$'s are negatively correlated (dependent only in that when a card moves in one attached tree, no card in another attached tree can move at the same time; in continuous time we would have had full independence)
\[\Var(S) = O(n^{7/3}\ep^{-6}\log^7n).\]
Thus $S \geq B_0n^{4/3}\ep^{-3}\log^3 n$ whp, by Chebyshev's inequality, for a suitable constant $B_0$.
However, each term $\tau_{r_i}$ of $S$ is, whp as noted above, the sum of at most $2\log n$ independent $C(h,l)$ conditioned on $\tau_l<\tau_r$ (plus an extra unconditioned $\tau_{r_i}$ at the end). Each such $C(h,l)$ is by the above arguments also bounded in distribution from above by $\E C(h^*,l)$ times an independent geometric$(1/2)$ number.
By the Chernoff bound, this means that $\Pro(\tau_{r_i} \leq 6\log n \E C(h^*,l)) = 1-o(1/n)$, i.e.\ whp no $\tau_{r_i}$ exceeds $6 \log n \E C(h^*,l) \leq B_1 n\ep^{-3}\log^3n$ for a constant $B_1$.
Hence at least $(B_0/B_1)n^{1/3}$ of the terms of $S$ must be at least $B_0 n \ep^{-3}\log^3 n$. This means that at time $B_0n\ep^{-3}\log^3n$ whp $\Theta(n^{1/3})$ cards $c_i$ will be in the same attached tree as they started in. At stationarity, the expected number of such cards is less than $1$. This proves that
\[\taum \geq B_0n\ep^{-3}\log^3 n\]
as desired.

\smallskip

\noindent {\bf Upper bound.}
Consider first the random multigraph, $H$, constructed in (1).
By a result of \cite{BFU}, on any bounded-degree graph $n_0$-vertex expander, one can choose
paths between each pair $(i,j)$, $1 \leq i < j \leq n_0$, of vertices.
such that no edge is used by more than $Bn_0\log n_0$ of these paths, where $B$ is some constant.
Now we want to apply this to $H$. It is well-known that $H$ is whp an expander. However, it is not bounded-degree. On the other hand, if we replace each vertex $v$ of $H$, with $d_v$ greater than, say, $D$, with a $\lceil d_v/D \rceil$-vertex expander and divide the edges to $v$ equally and uniformly at random among these neighbors (modulo integer considerations), we still whp  have an expander with $O(n)$ vertices.
In the sequel, we assume that $H$ has these low congestion paths and condition on that event.

Now extend $H$ to $K$ and then to $F$ according to (2) and (3) above.
Let each edge of $H$ be associated with one of its end vertices and let each vertex $v_1$ of $F$ be associated with the vertex of $H$, that the edge in $H$ from which $v_1$, or the attached tree containing $v_1$, arose, is associated with.
Let $X_v$, $v \in V(H)$, be the number of vertices of $F$ that are in this way associated with $v$.
Construct paths between each pair $(u_1,v_1)$ of vertices of $F$ by concatenating
the path chosen above between $u$ and $v$, with the shortest paths between $u_1$ and $u$ and $v_1$ and $v$, where $u$ and $v$ are the vertices of $H$ that $u_1$ and $v_1$ are associated with, respectively.
Let $L_e$ be the number of such paths that contain $e$, $e \in E(K)$.
Fix an $e \in E(K)$ and let $P_1,\ldots,P_A$, $A \leq BM\log M$, be an enumeration of the chosen paths of $H$ that contain $e$.
For convenience, define arbitrarily for each $P_i$, one of $P_i$'s end vertices as the start vertex and the other as the end vertex.
Let $u_i$ and $v_i$ be $P_i$'s start and end vertex respectively.
Then
\beq \label{ea}
L_e \leq \sum_{i=1}^A X_{u_i}X_{v_i}.
\eeq
Each $X_v$ is a sum of at most $C$ independent random variables, $Y_j(v)$,
where each $Y_j(v)$ is the sum of an independent geometric($\delta$) number of
independent random variables having the distribution of $T$, the total size
of a Poisson($1-\delta$) Galton-Watson tree.
Here $j$ ranges over $E(v)$, the set of
edges of $H$ associated with $v$.
Hence
\[L_e \leq \sum_{i=1}^A\sum_{j \in E(u_i)}\sum_{k \in E(v_j)}Y_j(u_i)Y_k(v_i).\]
By Cauchy-Schwarz,
\beq \label{ec}
L_e \leq \Big(\sum_{i=1}^A\sum_{j \in E(u_i)}\sum_{k \in E(v_i)}Y_j(u_i)^2\Big)^{1/2}\Big(\sum_{i=1}^A\sum_{j \in E(u_i)}\sum_{k \in E(v_i)}Y_k(v_i)^2\Big)^{1/2}.
\eeq
In each of the two sums on the right hand side of (\ref{ec}), each $Y_j(v)$ may appear with multiplicity, arising from that the corresponding $X_v$ may appear several times in (\ref{ea}).
It follows that the multiplicity is bounded by $2|E(H)|$ which in turn is bounded whp by $8M$.
Fix one of the right hand side sums arbitrarily.
Let $S_l$, $l=1,\ldots,8M$ be the set of $Y_j(v)$'s that
appears in the fixed sum with multiplicity at least $l$.
Then the $S_l$'s are decreasing and our fixed sum can be written in a more neat form as
\beq \label{ed}
\sum_{l=1}^{8M}\sum_{i \in S_l}Z_i^2
\eeq
where $Z_1,\ldots,Z_{8M}$ are iid with the distribution of $Y_j(v)$.
To bound the size of the $Z_j$'s we need to control the total size $T$ of an attached tree and sums of independent $T$'s.
The following is standard and can be found e.g.\ in \cite{Jagers}.

By the usual encoding of a Galton-Watson process into a random walk, it follows that $T$ has the distribution of the time of the first visit to $0$ of a simple random walk, starting from 1 and with step distribution Poisson$(1-\delta)-1$.
Hence the sum $Z=\sum_{r=1}^W T_r$, where the $T_r$'s are independent and distributed as $T$, has the distribution of the first visit to 0, starting from $W$.
The $Z_i$'s  in (\ref{ed}) have this form, with $W$ being
geometric($\delta$).
After $j\delta^{-2}$ steps of such a random walk, the position of the walker has
expectation $W-j\delta^{-1}$ and standard deviation less than
$\sqrt{j}\delta^{-1}$. By the Central Limit Theorem, the probability
that the walker at this time is above the $x$-axis, is bounded
by $e^{-(\sqrt{j}-\delta W/\sqrt{j})_{+}^2/2}$.
Hence
\[\Pro(Z \geq j\delta^{-2}|W) \leq e^{-(\sqrt{j}-\delta W/\sqrt{j})_{+}^2/2}
\leq e^{-(\sqrt{j}-\delta W/\sqrt{j})_{+}^2/4}
\leq e^{\delta W/2-j/4}.\]
Taking expectations gives
\[\Pro(Z \geq j\delta^{-2}) \leq \E[e^{\delta W/2}]e^{-j/4}.\]
The probability generating function of $W$ is $\delta/(1-s(1-\delta))$.
Plugging in $s=e^{\delta/2} = (1+o(1))(1+\delta/2)$, gives
\beq \label{ee}
\Pro(Z \geq j\delta^{-2}) \leq (1+o(1))2e^{-j/4} < e^{1-j/4}.
\eeq

Next we turn to the sum $\sum_l\sum_{i \in S_l}Z_i^2$.
Fix a large integer $s_0$.
For $|S_l| \leq s_0$, we have by (\ref{ee}) that the probability that $Z_j \geq s_0\delta^{-2}$ for any $j \in S_l$ is bounded by $s_0e^{1-s_0/4}$. Hence, for such $l$, we have
$\sum_{j \in S_l}Z_j^2 \leq s_0^3\delta^{-4}$ with probability at least $1-s_0e^{1-s_0/4}$.

Consider now $s := |S_l| \geq s_0+1$.
Order the $Z_i$'s in $S_l$: $Z_{(1)} \geq \ldots \geq Z_{(s)}$.
We have
\[\Pro(Z_{(r)} \geq j\delta^{-2}) \leq \Pro({\rm Bin}(s,e^{1-j/4}) \geq r)
\leq \frac{s^r e^{r-rj/4}}{r!}<s^{-3}\]
whenever $j \geq 4((1+3/r)\log s - \log r+2)$.
Hence $Z_{(r)} \leq 4\delta^{-2}((1+3/r)\log s - \log r+2)$ for all $r=1,\ldots,s$ with probability at least $1-s^{-2}$.
Hence, with probability at least $1-s^{-2}$, by standard integral approximation,
\beqa
\sum_{j \in S_l}Z_j^2 &\leq& 16\delta^{-4}\sum_{r=1}^s \Big((1+\frac3r)\log s - \log r + 2\Big)^2 < (160 s+O(\log^3 s))\delta^{-4} \\
&<& 161 s\delta^{-4}
\eeqa
provided that $s_0$ is large enough.
The probability that this holds for all $s=s_0+1,\ldots,8M$ is hence at least
$1-\sum_{s=s_0+1}^\infty s^{-2} > 1-\frac{1}{s_0}$.
Summing up in (\ref{ee}), using (\ref{ed}) and noting that there are whp $O(M\log M)$ terms in each of the sums on the right hand side of (\ref{ec}), we get
\beqa
L_e &\leq& \sum_{l:|S_l| \leq s_0}s_0^3\delta^{-4} + \sum_{l:|S_l| > s_0} 161|S_l|\delta^{-4} \leq \Big(s_0^4+161\sum_{l=1}^{8M}|S_l|\Big)\delta^{-4} \\
 &=& O(\delta^{-4}M\log M) = \Theta(\ep^{-2}n \log(\ep^2 n))
\eeqa
with probability at least $1-2(s_0e^{-s_0}+1/s_0)$, i.e.\ with probability arbitrarily close to $1$ on choosing $s_0$ sufficiently large.
This bounds the number of paths connecting all pairs of vertices of $F$ that
an edge of $K$ can appear in.
It remains to check in how many paths an edge of an attached tree
can appear.
However, given an attached tree, the edge that must be used most frequently is one which is incident to the vertex of $K$.
This edge is used in less than $tn$ paths, where $t$ is the size of the attached tree.
By the above, the largest tree has whp size $\Theta(\ep^{-2}\log(\ep^2 n))$, so no edge of an attached tree is used by more than $\Theta(\ep^{-2}n\log(\ep^2 n))$
times.
By the above, we also know that whp the longest path in $F$
is of order $\ep^{-1}\log(\ep^2n)$. Altogether, whp
\[A_* = \Theta(\ep^{-3}\log^2(\ep^2n)).\]
Plugging this into the comparison lemma gives whp
\[\taum \leq \tauh = O(n\ep^{-3}\log^2(\ep^2 n)\log n)\]
as desired.
\hfill $\Box$

\bigskip

\noindent {\bf Proof of Theorem \ref{tA}(e).}
This part is essentially a simplification or the proof of part (d).
It is easy to show that the giant component $G=(V,E)$ of $G(N,(1+\ep)/N)$ whp contains $\Omega(n^{2/3})$ sticks (i.e.\ induced paths with one end vertex having degree $1$) of length at least $\alpha\log n$ if the constant $\alpha$ is taken small enough.
See e.g.\ \cite{FR}.
Letting these sticks play the role of the attached trees $T_i$, in the lower bound part in the proof of (d), establishes that
\[\taum = \Omega(n\log^3n).\]
Now turn to the upper bound.
The construction (1)-(3) above does not in this case exactly generate a graph contiguous to $G$. There is however, a characterization of $G$, due to Benjamini et al.\ \cite{BKW} that
works equally well for our purposes. According to Theorem 4.2 of \cite{BKW}, there is an $\alpha \in (0,1)$ such that whp $G$ contains a subgraph $B$ such that
\bi
\item[(i)] $B$ is an $\alpha$-expander (i.e.\ for any $W \subset V$ with $|W| \leq n/2$, the number of edges connecting $W$ to $V$ is at least $\alpha|W|$).
\item[(ii)] Letting $D_i$ denote the connected components of $G \setminus B$ and $E'(D_i)$ the set of edges with at least one endpoint in $D_i$, the number of $i$ for which $|E'(D_i)| \geq k$ is bounded above by $n(1+\ep)e^{-\alpha k}$, $k=1,2,\ldots$.
\item[(iii)] Each $v \in B$ is connected to no more than $1/\alpha$ different $D_i$'s
\ei
Using (i), we have as above that all vertex pairs of $B$ can be joined by paths in $B$ in such a way that no edge of $B$ is used in more than $O(n\log n)$ paths.
Then, by (ii) and (iii), extending this to joining all pairs of $G$ with paths such that no edge of $G$ is used more than $O(n \log n)$ times, goes through in a much simplified form,
without most of the calculations involved.
Plugging this into the comparison technique is done exactly as before and yields
\[\taum = O(n\log^3n).\]
\hfill $\Box$

\bigskip

\noindent {\bf Proof of Theorem \ref{tA}(f).} Here an appeal to Border et al.\ \cite{BFU}, establishing that all pairs of vertices can be joined by paths so that no edge is used more than
$O(\log n)$ times, goes directly into the comparison technique, giving $\taum = O(n\log^3 n)$.

The lower bound follows exactly as the standard argument for the random transpositions shuffle: after time time $(1/2-a)n\log n$ whp at least $n^{a/2}$ cards have never bee moved and are hence in their starting positions. However at stationarity, the expected number of such cards is $1$, so the stationary probability for this event is $o(1)$.
\hfill $\Box$

\bigskip

\noindent {\bf Proof of Theorem \ref{tA}(g).} The lower bound was proved by Wilson \cite{Wilson1}. Indeed, with $\xi(v)=\xi(v_1,\ldots,v_d) = 2v_1-1$, $\xi$ becomes an eigenvector to an eigenvalue of order $1/(n\log n)$ for the walk of a single card. Thus Proposition \ref{pB} establishes the lower bound.

For the upper bound, we again use the comparison technique, comparing with the random transpositions shuffle.
Since there are $\Theta(n^2)$ pairs of vertices, the diameter of $\Z_2^d$ is $d$ and the number of edges in $\Z_2^d$ is $\Theta(nd)$, any automorphism invariant rule for assigning representations of arbitrary transpositions via shortest paths will have that any edge appears in $\Theta(n)$ of these paths.
Now use the comparison technique as before.
\hfill $\Box$

\end{document}